\documentclass[11pt]{amsart}

\usepackage{amsfonts, amstext, amsmath, amsthm, amscd, amssymb}
\usepackage{epsfig, graphics, color}
\usepackage{multicol}

\usepackage[margin=1.2in]{geometry}

\let\oldmarginpar\marginpar
\renewcommand\marginpar[1]{\oldmarginpar[\raggedleft\footnotesize #1]%
{\raggedright\footnotesize #1}}

\renewcommand{\setminus}{{\smallsetminus}}

\newcommand{\vol}{{\rm vol}}

\theoremstyle{plain}
\newtheorem{theorem}{Theorem}[section]

\newtheorem{lemma}[theorem]{Lemma}

\newtheorem*{namedtheorem}{\theoremname}
\newcommand{\theoremname}{testing}
\newenvironment{named}[1]{\renewcommand{\theoremname}{#1}\begin{namedtheorem}}{\end{namedtheorem}}

\theoremstyle{definition}

\newtheorem{remark}[theorem]{Remark}

\begin{document}
\title{Volumes of chain links}
\author[J. Kaiser]{James Kaiser}
\author[J. Purcell]{Jessica S. Purcell}
\author[C. Rollins]{Clint Rollins}

\address[]{ Department of Mathematics, Brigham Young University,
Provo, UT 84602}



\begin{abstract}
Agol has conjectured that minimally twisted $n$--chain links are the
smallest volume hyperbolic manifolds with $n$ cusps, for $n \leq 10$.
In his thesis, Venzke mentions that these cannot be smallest volume
for $n \geq 11$, but does not provide a proof.  In this paper, we give
a proof of Venzke's statement for a number of cases. For $n \geq 60$
we use a formula from work of Futer, Kalfagianni, and Purcell to
obtain a lower bound for volume. The proof for $n$ between $12$ and
$25$ inclusive uses a rigorous computer computation that follows
methods of Moser and Milley.  Finally, we prove that the $n$--chain
link with $2m$ or $2m+1$ half--twists cannot be the minimal volume
hyperbolic manifold with $n$ cusps, provided $n \geq 60$ or $|m| \geq
8$, and we give computational data indicating this remains true for
smaller $n$ and $|m|$.
\end{abstract}

\maketitle

\section{Introduction}\label{sec:intro}
\begin{figure}[h!]
  \includegraphics{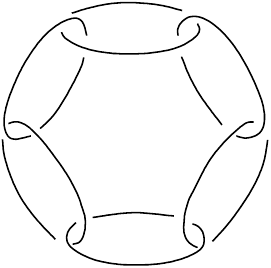}
  \hspace{.5in}
  \includegraphics{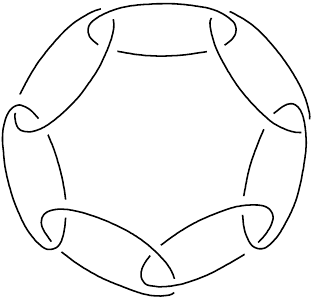}
  \hspace{.5in}
  \includegraphics{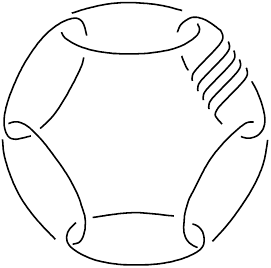}
  \caption{Left: Minimally twisted $6$--chain link.  Middle:
    Minimally twisted $7$--chain link.  Right: $6$--chain
    link with more half--twists.}
  \label{fig:chain-link1}
\end{figure}

An \emph{$n$--chain link} consists of $n$ unknotted circles embedded
in $S^3$, linked together in a closed chain.  Notice that links of a
chain can be connected with an arbitrary amount of twisting.  In
particular, if we embed the first link in the plane of projection, the
next perpendicular to the plane of projection, the next again in the
plane of projection, and so on, then the last link may include any
integer number of half--twists.  See, for example, Figure
\ref{fig:chain-link1}.

Hyperbolic structures on $n$--chain link complements have been
studied, for example, by Neumann and Reid
\cite{neumann-reid:arithmetic}.  They show any $n$--chain link
complement with $n\geq 5$ admits a hyperbolic structure. In this
paper, we are primarily interested in hyperbolic manifolds, so we
restrict our attention to $n\geq 5$.

A \emph{minimally twisted $n$--chain link} is an $n$--chain link such
that, if $n$ is even, each link component alternates between lying
embedded in the projection plane and lying perpendicular to the
projection plane.  If $n$ is odd, the link may be arranged such that
each component alternates between lying in the projection plane and
perpendicular to it, except a single link component which connects a
link which is embedded in the projection plane to one which is
perpendicular, with no twisting.  See Figure \ref{fig:chain-link1}.

Notice that there are actually two choices for the minimally twisted
$n$--chain link for $n$ odd, depending on which way the last links are
connected.  However, these are isometric by an orientation reversing
isometry, so we will not distinguish between them. 

In \cite{agol:min-vol}, Agol conjectures that minimally twisted
$n$--chain link complements, for $n\leq 10$, are the smallest volume
hyperbolic 3--manifolds with exactly $n$ cusps, but notes that Venzke
has pointed out they cannot be smallest for $n\geq 11$, as the
$(n-1)$--fold cyclic cover over one component of the Whitehead link
has smaller volume.  This statement is included in Venzke's thesis
\cite{venzke:thesis}.  However, Venzke does not give a proof.  In this
paper, we give a rigorous proof for $n\geq 60$.  The following theorem
is the main result of this paper.

\begin{named}{Theorem \ref{thm:not-minvolume}}
For $n\geq 60$, the minimally twisted $n$--chain link complement has
volume strictly greater than that of the $(n-1)$--fold cyclic cover
over one component of the Whitehead link.  Hence the minimally twisted
$n$--chain link complement cannot be the smallest volume hyperbolic
3--manifold with $n$ cusps, $n\geq 60$.
\end{named}

For $n$ between $11$ and $59$, inclusive, we present computer
tabulation of volumes, compared with the volumes of the $(n-1)$--fold
cyclic cover of the Whitehead link.  See Table \ref{table:volumes}.
By inspection, Theorem \ref{thm:not-minvolume} also holds for these
manifolds.  In Section \ref{sec:compute}, we explain how these
computations can be made completely rigorous --- at least for those
values of $n$ for which the minimally twisted $n$--chain link
complement is triangulated with fewer than $100$ tetrahedra.  In this
case, this includes $n$ between 12 and 25, inclusive.

Finally, in Section \ref{sec:arblinks} we evaluate volumes of
arbitrarily twisted $n$--chain link complements.  The main result of
that section is Theorem \ref{thm:vol-withtwist}, which states that no
$n$--chain link complement can be the minimal volume $n$--cusped
hyperbolic 3--manifold, provided either $n\geq 60$, or the chain link
contains at least $17$ half--twists.  We present computational data to
show that similarly, for $11 \leq n \leq 59$, no $n$--chain link
complement can be minimal volume.  When $5 \leq n \leq 10$, we
rigorously prove, by computer, that no $n$--chain link which is
\emph{not} minimally twisted can be the minimal volume $n$--cusped
hyperbolic 3--manifold.

\begin{remark}
Since a version of this paper was made public, Hidetoshi Masai has
pointed out to us that for even chain links, more can be said.  In
\cite[Chapter 6]{thurston}, Thurston finds a formula for the volumes
of minimally twisted chain links with an even number of link
components. Namely,
$$\vol(S^3\setminus C_{2n}) = 8n \left( \Lambda\left(\frac{\pi}{4} +
\frac{\pi}{2n}\right) + \Lambda\left(\frac{\pi}{4} -
\frac{\pi}{2n}\right)\right),$$ where $\Lambda$ is the Lobachevsky
function.  Masai notes that the difference of this volume and the
volume of the $(2n-1)$--cyclic cover over a component of the Whitehead
link is an increasing function in $n$, for $n \geq 6$.  This result
will give a rigorous proof that minimally twisted $n$--chain links are
not minimal volume for $17$ additional chain links, namely for $n =
26, 28, 30, \dots, 56, 58$.  In addition, this gives an alternate
proof that minimally twisted $2n$--chain links are not minimal volume
for larger $n$.  However, the result for odd links requires other
techniques, for instance those in this paper.
\end{remark}  

\subsection{Acknowledgements}
We thank Ian Agol and Hidetoshi Masai for helpful conversations.  We
also thank Peter Milley, Harriet Moser, and Saul Schleimer for their
assistance with computational aspects of this project.  Authors Kaiser
and Rollins were both supported in part by Brigham Young University
mentoring funds for undergraduate research.  All authors were
supported in part by a grant from the National Science Foundation.

\section{Slope lengths on covers}\label{sec:slopes}

To prove the main theorems of this paper, for $n\geq 60$, we will
obtain the complement of the minimally twisted $n$--chain link by Dehn
filling a manifold $\widehat{W}_n$ which is geometrically explicit,
constructed by gluing together manifolds isometric to the Whitehead
link complement, cut along 2--punctured disks.

\begin{figure}
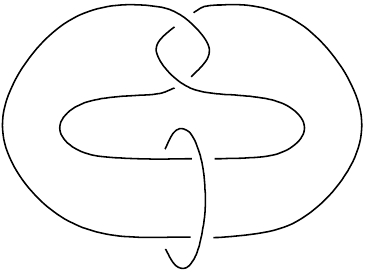 \hspace{.5in}
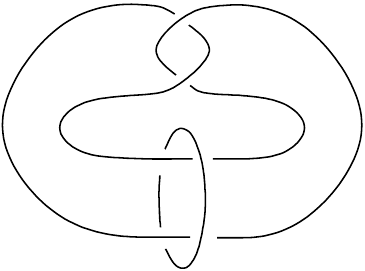
\caption{The Whitehead link (left), and its reflection (right), with
component labeled $K$ and $\overline{K}$, respectively.}
\label{fig:w1}
\end{figure}

We work with the diagram of the Whitehead link as in Figure
\ref{fig:w1}, left, with a link component denoted $K$.  Note that by
switching the direction of a pair of crossings, we obtain a link whose
complement is isometric to that of the Whitehead link complement by an
orientation reversing isometry.  The isometry takes $K$ to a link
component we denote $\overline{K}$, as in Figure \ref{fig:w1}, right.

\begin{lemma}
The shape of the cusp of $K$ is a parallelogram with meridian and
longitude meeting at angle $-\pi/4$ (measured from meridian to
longitude).

Similarly, the shape of the cusp of $\overline{K}$ is a parallelogram with
meridian and longitude meeting at angle $\pi/4$ (measured from meridian
to longitude).

When we take a maximal horocusp about $K$ or $\overline{K}$, the
meridian has length $\sqrt{2}$, and the longitude has length $4$.
\label{lemma:k_shape}
\end{lemma}

\begin{figure}
\includegraphics{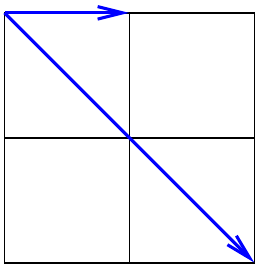}
\caption{Cusp shape of component $K$ of Whitehead link.  Meridian runs
horizontally across the top, longitude runs diagonally.}
\label{fig:whitehead-cusp}
\end{figure}

Lemma \ref{lemma:k_shape} is illustrated in
Figure \ref{fig:whitehead-cusp}.

\begin{proof}
The first part of the lemma, and lengths of slopes on $K$, are well
known for the Whitehead link.  See, for example
\cite{neumann-reid:arithmetic}.

As for $\overline{K}$, the orientation reversing isometry taking the
Whitehead link to its reflection takes a meridian of $K$ to a meridian
of $\overline{K}$, and reflects the longitude.  Since this is an
isometry, the meridian and longitude of $\overline{K}$ have the same
lengths as those of $K$, but the angle from the meridian to longitude
is reflected across the meridian, to be $\pi/4$.
\end{proof}

Now, the manifold $\widehat{W}_n$ can be described as the minimally
twisted $n$--chain link embedded in a standard solid torus.
Therefore, to obtain the complement of the minimally twisted
$n$--chain link, we will Dehn fill $\widehat{W}_n$ along a standard
longitude of the solid torus boundary component.

\begin{figure}
\includegraphics{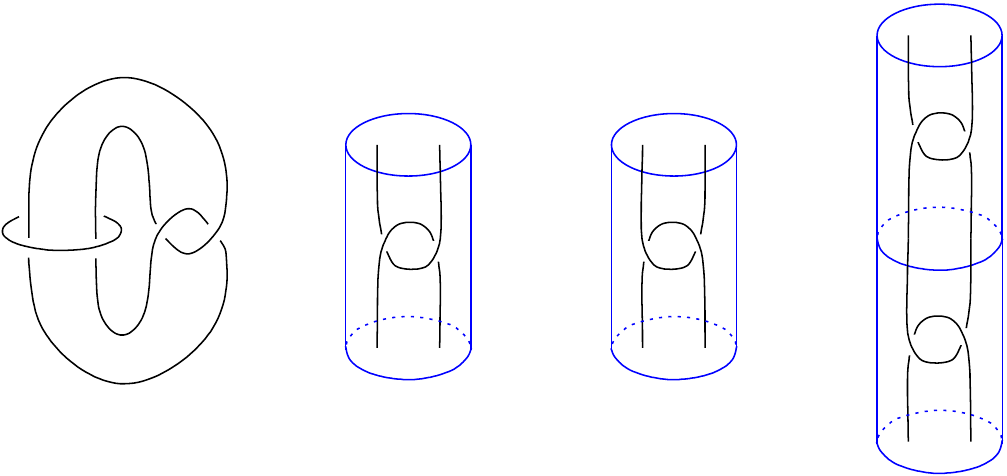}
\caption{Constructing the manifold $\widehat{W}_n$: cut the Whitehead
link complement along a 2--punctured disk to obtain $W_1$ (second from
left).  Its reflection is $\overline{W}_1$ (third from left).  Attach
these to form a link in a solid cylinder, right.}
\label{fig:wn-construct}
\end{figure}

To build $\widehat{W}_n$ from the Whitehead link complement, proceed
as follows.  First, cut the Whitehead link complement along the
2--punctured disk bounded by $K$ to get a clasp in a cylinder, which
we call $W_1$, shown second from left in Figure
\ref{fig:wn-construct}.  Similarly, cut the reflected Whitehead link
complement along the 2--punctured disk bounded by $\overline{K}$ to
get a clasp in the opposite direction in a cylinder, which we call
$\overline{W}_1$, shown third from left in Figure
\ref{fig:wn-construct}.

Now, attach a copy of $W_1$ to $\overline{W}_1$ via an isometry of the
2--punctured disk as on the right in Figure \ref{fig:wn-construct}.
In particular, boundary components are glued as shown without
twisting.  Call the resulting link in a solid cylinder $W$.

For $n$ even, glue ${n}/{2}$ copies of $W$ together end to end,
without twisting, followed by gluing the remaining two ends.  For $n$
odd, glue $(n-1)/{2}$ copies of $W$ together without twisting,
then glue a single copy of $W_1$, and attach the ends without
twisting.  This completes the construction of $\widehat{W}_n$.

\begin{lemma}
Let $\epsilon = n\mod{2}$.  The minimally twisted $n$--chain link in a
solid torus, $\widehat{W}_n$, has solid torus boundary component
comprised of $\lfloor n/2 \rfloor + \epsilon$ copies of the
cusp $K$ coming from $W_1$, and $\lfloor n/2 \rfloor$ copies
of the cusp $\overline{K}$ coming from $\overline{W}_1$.  The standard
longitude of the solid torus follows a meridian of each copy of $K$
and $\overline{K}$, where the meridians of each copy of $K$ are
orthogonal to the meridians of each copy of $\overline{K}$.  The
length of the longitude of the solid torus boundary component is
$\sqrt{n^2 + \epsilon}$.
\label{lemma:torus_boundary_component}
\end{lemma}

\begin{proof}
This follows from the construction of $\overline{W}_n$ and
Lemma \ref{lemma:k_shape}.

The boundary component corresponding to the solid torus comes from
$n/2$ copies of the cusp $K$ and $n/2$ copies of the
cusp $\overline{K}$ for $n$ even, and $(n-1)/2 + 1$ copies of
the cusp $K$ and $(n-1)/2$ copies of the cusp $\overline{K}$,
for $n$ odd.  These are glued together along their respective
longitudes.  Since a copy of the cusp of $K$ is glued to one of
$\overline{K}$ along the longitude of each, the meridians meet at
right angles.  See Figure \ref{fig:cusp-wnbar}.

\begin{figure}
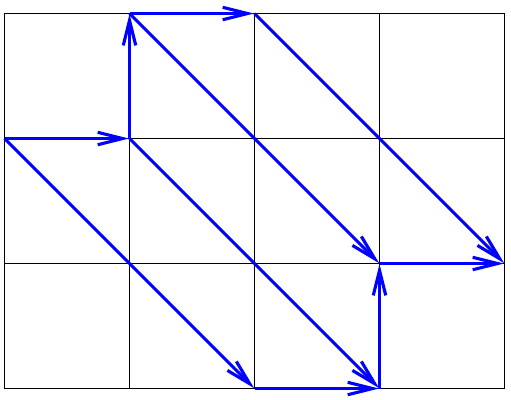
  \caption{The cusp construction of $\widehat{W_n}$.}
  \label{fig:cusp-wnbar}
\end{figure}

The longitude of the solid torus of $\widehat{W_n}$ is given by
following each meridian of the copies of $K$ and $\overline{K}$ that
glue to give the solid torus boundary component.  Since these
meridians always meet at right angles, the length of the longitude of
the solid torus may be determined by the Pythagorean theorem.  By
Lemma \ref{lemma:k_shape}, the length of the meridian of the cusp $K$,
and that of the cusp $\overline{K}$, is $\sqrt{2}$. We see that the
length of the longitude of the solid torus boundary component of
$\widehat{W}_n$ is $\sqrt{\left(\sqrt{2} \left( \lfloor \frac{n}{2} \rfloor
  + \epsilon \right)\right)^2 + \left( \sqrt{2} \lfloor \frac{n}{2} \rfloor
  \right)^2} = \sqrt{n^2 + \epsilon}$.
\end{proof}

Notice that while the construction of $\widehat{W}_n$ as described
above uses $\lfloor n/2 \rfloor + \epsilon$ copies of $W_1$
and $\lfloor n/2 \rfloor$ copies of $\overline{W}_1$, we could
have constructed it using $\lfloor n/2 \rfloor$ copies of
$W_1$ and $\lfloor n/2 \rfloor + \epsilon$ copies of
$\overline{W}_1$ instead.  The result of this modified method of
construction would be isometric to that of the original construction,
by an orientation reversing isometry.

To obtain $C_n$ from $\widehat{W}_n$ we Dehn fill a slope on the solid
torus boundary component of $\widehat{W}_n$ that follows one standard
longitude of the solid torus.  Hence we Dehn fill along a slope of
length $\sqrt{n^2 + (n\mod{2})}$.

\section{Volumes, large minimally twisted chains}\label{sec:volumes}

Using the information on slopes above, we may deduce geometric
information on minimally twisted chain links by applying appropriate
theorems bounding change in geometry under Dehn filling.  This will
give the desired result when $n\geq 60$.

\subsection{Dehn filling and volume}

We use the following theorem, which is a slightly simpler version of
the main theorem in \cite{fkp:dfvjp}.

\begin{theorem}[Futer--Kalfagianni--Purcell \cite{fkp:dfvjp}]
  Let $M$ be a complete, finite--volume hyperbolic manifold with (at
  least one) cusp, with horoball neighborhood $C$ about that cusp, and
  let $s$ be a slope on $\partial C$ with length $\ell(s) > 2\pi$.
  Then the manifold $M(s)$ obtained by Dehn filling $M$ along $s$ is
  hyperbolic, with volume:
  $$ \vol(M(s)) \: \geq \:
  \left(1-\left(\frac{2\pi}{\ell(s)}\right)^2\right)^{3/2} \vol(M).$$
\label{thm:fkp}
\end{theorem}

Putting this theorem together with Lemma
\ref{lemma:torus_boundary_component}, we obtain the following.

\begin{theorem}
For $n\geq 7$, the volume of the complement of the minimally twisted
$n$--chain link $C_n$ is at least:
$$\vol(S^3\setminus C_n) \geq n\, v_8 \, \left( 1 - \frac{4\pi^2}{ n^2
  + \epsilon } \right)^{3/2},$$
where $v_8 = 3.66386\ldots$ is the volume of a hyperbolic regular
ideal octahedron, and $\epsilon = n\mod{2}$.
\label{thm:volume}
\end{theorem}

\begin{proof}
The chain link complement is obtained by Dehn filling the longitude of
$\widehat{W_n}$.  This is obtained by gluing copies of the Whitehead
link and its reflection along totally geodesic 3--punctured spheres,
hence has volume $n$ times the volume of the Whitehead link, $n\cdot v_8$.

By Lemma \ref{lemma:torus_boundary_component}, we know the length of
the Dehn filling slope is $\sqrt{n^2 + \epsilon}$.  The result follows
by putting this data into Theorem \ref{thm:fkp}.
\end{proof}

We now give a proof of the main theorem.  

\begin{theorem}
For $n\geq 60$, the minimally twisted $n$--chain link complement has
volume strictly greater than that of the $(n-1)$--fold cyclic cover
over one component of the Whitehead link.  Hence the minimally twisted
$n$--chain link complement cannot be the smallest volume hyperbolic
3--manifold with $n$ cusps, $n\geq 60$.
\label{thm:not-minvolume}
\end{theorem}

\begin{proof}
The volume of the $(n-1)$--fold cyclic cover of the Whitehead link is
$(n-1)\,v_8$.  By Theorem \ref{thm:volume}, the volume of the
complement of the minimally twisted $n$--chain link is
$$\vol(S^3\setminus C_n) \geq
n\,v_8\,\left(1-\frac{4\pi^2}{n^2+\epsilon}\right)^{3/2} \geq 
n\,v_8\,\left(1-\frac{4\pi^2}{n^2}\right)^{3/2}.$$
We want to find $n$ for which the following inequality holds:
$$n\,v_8\,\left(1-\frac{4\pi^2}{n^2}\right)^{3/2} > (n-1)\,v_8,$$
or
\begin{equation}\label{ineq}
  \left( \frac{n}{n-1} \right)\left(1-\frac{4\pi^2}{n^2}\right)^{3/2}
-1 > 0.
\end{equation}
Let $f(n)$ be the function on the left side of inequality
\eqref{ineq}.  Using calculus, one sees that
$\lim_{n\to\infty}f(n)=0$, $f$ is increasing between $n=7$ and
$n=6\pi^2 + 2\pi\sqrt{9\pi^2 -2} \approx 117.8$, and decreasing for
larger $n$, which implies $f$ has at most one root for $n\geq 7$, and
that $f$ is positive to the right of any root.  The Intermediate Value
Theorem implies that there is a root between $n=59$ and $n=59.1$.
Hence the inequality is satisfied for $n\geq 60$.
\end{proof}

\section{Computations of volume, smaller minimally twisted chains}\label{sec:compute}

Now we analyze volumes of minimally twisted $n$--chain links
for $n$ between $11$ and $59$, since the main method of proof of
Theorem \ref{thm:not-minvolume} will not apply to these manifolds.

For $n$ between $11$ and $59$ inclusive, in Table
\ref{table:volumes} we present computational data using SnapPea
(SnapPy) \cite{weeks:snappea, snappy} that shows that the minimally
twisted $n$--chain link complement cannot be the minimal volume
hyperbolic 3--manifold with $n$ cusps.  In particular, $W_{n-1}$, the
$(n-1)$--fold cyclic cover over one component of the Whitehead link,
has smaller volume.  The volume of $W_{n-1}$ is always $(n-1)\,v_8$,
where $v_8 = 3.66386\ldots$ is the volume of a hyperbolic regular
ideal octahedron, which is the volume of the Whitehead link
complement.  Notice that for $n\geq 11$, the volume of $S^3\setminus
C_n$ is strictly larger than that of $W_{n-1}$.

\begin{table}
  \begin{multicols}{2}

{\footnotesize
  
\begin{tabular}{ | c || c | c |}
  \hline
$n$ & $\vol(S^3\setminus C_n)$ & $\vol(W_{n-1})$ \\
  \hline
5	&	10.14941606	&	14.65544951	\\
6	&	14.65544951	&	18.31931188	\\
7	&	19.79685462	&	21.98317426	\\
8	&	24.09218408	&	25.64703664	\\
9	&	28.47566906	&	29.31089901	\\
10	&	32.55154031	&	32.97476139	\\
11	&	36.64918655	&	36.63862377	\\
12	&	40.59766426	&	40.30248614	\\
13	&	44.5536682	&	43.96634852	\\
14	&	48.42519197	&	47.63021090	\\
15	&	52.29990219	&	51.29407327	\\
16	&	56.12184477	&	54.95793565	\\
17	&	59.94533184	&	58.62179803	\\
18	&	63.73354269	&	62.28566041	\\
19	&	67.52257845	&	65.94952278	\\
20	&	71.28681886	&	69.61338516	\\
21	&	75.05153335	&	73.27724753	\\
22	&	78.79813245	&	76.94110991	\\
23	&	82.54502011	&	80.60497229	\\
24	&	86.27825885	&	84.26883466	\\
25	&	90.01168157	&	87.93269704	\\
26	&	93.73455871	&	91.59655942	\\
27	&	97.45755771	&	95.26042179	\\
28	&	101.1722364	&	98.92428417	\\
29	&	104.8869984	&	102.5881465	\\
30	&	108.5950782	&	106.2520089	\\
31	&	112.3032167	&	109.9158713	\\
32	&	116.0059062	&	113.5797337	\\
\hline
\end{tabular}

\begin{tabular}{|c||c|c|}
  \hline
$n$ & $\vol(S^3\setminus C_n)$ & $\vol(W_{n-1})$ \\
  \hline
33	&	119.708638	&	117.2435961	\\
34	&	123.4068675	&	120.9074584	\\
35	&	127.1051279	&	124.5713208	\\
36	&	130.7996249	&	128.2351832	\\
37	&	134.494145	&	131.8990456	\\
38	&	138.1854868	&	135.5629079	\\
39	&	141.8768462	&	139.2267703	\\
40	&	145.5654969	&	142.8906327	\\
41	&	149.2541611	&	146.5544951	\\
42	&	152.9404979	&	150.2183574	\\
43	&	156.6268453	&	153.8822198	\\
44	&	160.3111779	&	157.5460822	\\
45	&	163.995519	&	161.2099446	\\
46	&	167.6781044	&	164.8738070	\\
47	&	171.3606965	&	168.5376693	\\
48	&	175.0417493	&	172.2015317	\\
49	&	178.7228075	&	175.8653941	\\
50	&	182.4025087	&	179.5292565	\\
51	&	186.0822143	&	183.1931188	\\
52	&	189.7607173	&	186.8569812	\\
53	&	193.4392239	&	190.5208436	\\
54	&	197.1166599	&	194.1847060	\\
55	&	200.7940988	&	197.8485683	\\
56	&	204.4705802	&	201.5124307	\\
57	&	208.1470642	&	205.1762931	\\
58	&	211.8226885	&	208.8401555	\\
59	&	215.4983149	&	212.5040178	\\
60	&	219.1731666	&	216.1678802	\\
\hline
\end{tabular}
}					
					
\end{multicols}

  \caption{Volumes of the complement of the minimally twisted
    $n$--chain link $C_n$, compared to volumes of $W_{n-1}$, the
    $(n-1)$--fold cyclic cover over a component of the Whitehead link,
    for $5 \leq n \leq 60$.  Note that $S^3\setminus C_n$ has greater
    volume for $n\geq 11$.}
  \label{table:volumes}
\end{table}

It would be nice to turn this data into a rigorous proof that the
minimally twisted $n$--chain links for $n$ between $11$ and $59$
cannot be minimal volume.  One way to do this would be to use the
methods of Moser \cite{moser} and Milley in \cite{milley:minvol}.
Milley has written a program to rigorously prove that a hyperbolic
3--manifold with hyperbolic structure computed by Snap
\cite{goodman:snap} has volume greater than some constant.  This
program, which is available as supplementary material with
\cite{milley:minvol}, is in theory exactly what we need for these
chain link examples.

However, in practice, making Moser and Milley's programs work with the
chain links has proven to be difficult, due to the computational
complexity of the chain links.  While Milley worked with small
manifolds, for example with less than $10$ tetrahedra, and Moser's
largest manifold included $57$ tetrahedra, our triangulations of
minimally twisted $n$--chain link complements include between $40$ and
$236$ tetrahedra.  We were successfully able to run Moser's algorithm
for $n$ between $11$ and $25$, inclusive, which gives results for
those manifolds triangulated with up to $100$ tetrahedra, but then the
program failed.  We were able to run Milley's algorithm for all values
of $n$ for which Moser's algorithm applied. However, Milley's
algorithm only returned a positive result for $n$ between 12 and 25,
inclusive.  Therefore, we have the following result.

\begin{theorem}\label{thm:vol-computational}
  For $n$ between $12$ and $25$, inclusive, the minimally twisted
  $n$--chain link complement has volume strictly greater than that of
  the $(n-1)$--fold cyclic cover over one component of the Whitehead
  link, hence cannot be the smallest volume hyperbolic 3--manifold
  with $n$ cusps.
\end{theorem}

\begin{proof}
The proof is identical to that of Milley \cite{milley:minvol}, and
uses his code, included as supplementary material with that reference
\cite{milley:minvol}, modified to read in minimally twisted $n$--chain
links rather than Dehn fillings of census manifolds.  The first step
is to feed the triangulations of the minimally twisted $n$--chain
links into Snap, and ensure that the triangulations used in the
computation are geometric, that is, all tetrahedra are positively
oriented.  This is true for all minimally twisted $n$--chain links,
$11 \leq n \leq 59$.

Next, use Moser's algorithm \cite{moser} to find a value $\delta$
which measures the maximal error between Snap's computed solution and
the true solution.  Moser's algorithm gave us such a value for $11
\leq n \leq 25$, but failed thereafter, presumably due to
computational complexity of the chain link complements.

Finally, for each $n$ between $12$ and $25$, inclusive, input the Snap
triangulation data and Moser's value $\delta$ into Milley's program
{\tt{rigorous\_volume.C}}, along with the constant value
$(n-1)*3.66386237670888$.  The program checks rigorously whether the
volume of the given $n$--chain link is larger than the given constant.
For $12 \leq n \leq 25$, the program definitively proved that the
volumes of the minimally twisted $n$--chain link complement are
strictly larger than that of the $(n-1)$--fold cyclic cover over a
component of the Whitehead link.
\end{proof}

\begin{remark}
Note that the above theorem does not hold for $n=11$.  Although the
triangulation of the minimally twisted 11--chain link is positively
oriented, and Moser's algorithm returns a value of $\delta$ for this
link, Milley's program {\tt{rigorous\_volume.C}} is unable to verify
that its volume is larger than that of the 10--fold cyclic cover of
the Whitehead link.  When $n=11$, the volumes of these manifolds are
too close for rigorous checking.
\end{remark}

What about the volumes output by SnapPea for $26 \leq n \leq 59$?
Note in Table \ref{table:volumes} that the minimally twisted
$n$--chain link for these $n$ has volume greater than $2$ plus the
volume of $W_{n-1}$.  It is highly unlikely that SnapPea's computation
would be so far off as to make the theorem untrue for any of these
values of $n$.  However, since we do not have a rigorous proof at this
time, we do not include the result as a theorem.

\section{Arbitrary chain links}\label{sec:arblinks}
In this section, we extend our results to chain links with an
arbitrary amount of twisting.

Consider again the manifold $\widehat{W}_n$, which is a minimally
twisted $n$--chain link in a solid torus.  Let $\lambda_n$ denote the
standard longitude of the solid torus, and let $\mu_n$ denote the
meridian.  In the previous section, we performed Dehn filling along
the slope $\lambda_n$ to obtain the complement of the minimally
twisted $n$--chain link in $S^3$.  Notice that Dehn filling along any
slope of the form $\lambda_n + m\, \mu_n$ will also yield the
complement of a chain link in $S^3$, where the resultant chain has
$|m|$ additional full twists (or $2|m|$ additional crossings).  The
twisting will be positive or negative depending on the sign of $m$.

\begin{lemma}\label{lemma:slope-length-nohalftwist}
  For any integer $m$, the slope $\lambda_n + m\, \mu_n$ on the solid
  torus boundary component of $\widehat{W}_n$ has length
  $\sqrt{n^2 + 16m^2 + (n\mod 2)(1+8m)}$.
\end{lemma}

\begin{proof}
The solid torus boundary component of $\widehat{W_n}$ is tiled by
regular ideal octahedra coming from the Whitehead link, and these
appear as squares of side length $\sqrt{2}$ by
Lemma \ref{lemma:k_shape}.  When we place the corner of one such
square at $(0,0)$ in the Euclidean plane, we see that the slope
$\lambda_n$ runs from $(0,0)$ to $(\sqrt{2}(\lfloor n/2 \rfloor +
(n\mod 2)), \sqrt{2}\lfloor n/2
\rfloor)$, as in Lemma \ref{lemma:torus_boundary_component}.

The slope $\mu_n$ runs along exactly one of the 2--punctured disks we
sliced in the Whitehead link complement (or its reflection) to build
$\widehat{W_n}$.  Hence by Lemma \ref{lemma:k_shape} it runs from
$(0,0)$ to $(2\sqrt{2}, -2\sqrt{2})$.

Thus the slope $\lambda_n + m\mu_n$ runs from $(0,0)$ to
$(\sqrt{2}(\lfloor n/2 \rfloor + (n\mod 2) + 2m, \sqrt{2}(\lfloor n/2
\rfloor - 2m)$, hence has length as follows.

For $n$ even, $n=2k$,
$$ \ell(\lambda_n + m\mu_n) = \left( 2(k + 2m)^2 + 2(k-2m)^2
\right)^{1/2} = \sqrt{n^2 + 16m^2}.
$$
For $n$ odd, $n=2k+1$,
$$ \ell(\lambda_n + m\mu_n) = \left( 2(k+1+2m)^2 + 2(k-2m)^2 \right) =
\sqrt{n^2 + 16m^2 + (1+8m)}.
$$
\end{proof}

In order to obtain any $n$--chain link, in addition to considering
Dehn filling on the manifold $\widehat{W}_n$, we must also consider
Dehn filling on an $n$--chain link in the solid torus which differs
from $\widehat{W}_n$ by the insertion of a single crossing, or
half--twist, at a 2--punctured disk.  We call this manifold
$\overline{W}_n$.  Recall that we constructed $\widehat{W}_n$ by
gluing together alternating copies of $W_1$ and $\overline{W}_1$ along
their 2--punctured disk boundaries, without twisting, to form a link
in a solid cylinder, and then gluing the cylinder end to end without
twisting.  To form $\overline{W}_n$, we may glue by a single
half--twist when we connect the final solid cylinder end to end.
Equivalently, if $n$ is even, replace the last copy of
$\overline{W}_1$ with $W_1$ and glue end to end without twisting.  If
$n$ is odd, replace the last $W_1$ with $\overline{W}_1$ and glue end
to end without twisting.  This gives the desired half--twist in both
cases.

Now, denote the standard longitude of the solid torus boundary
component of $\overline{W}_n$ by $\bar{\lambda}_n$, and the meridian
by $\bar{\mu}_n$.  Dehn filling along a slope of the form
$\bar{\lambda}_n + m\,\bar{\mu}_n$ yields the complement of an
$n$--chain link in $S^3$, which differs from the minimally twisted
$n$--chain link by the insertion of $2m + 1$ half--twists, where the
direction of half--twist is determined by the sign of $m$.

\begin{lemma}\label{lemma:slope-length-halftwist}
  For any integer $m$, the slope $\bar{\lambda}_n + m\,\bar{\mu}_n$ on
  the solid torus boundary component of $\widehat{W}_n$ has length
  $\sqrt{n^2 + 4(1+2m)^2}$, if $n$ is even, and $\sqrt{n^2 + 16m^2 +
    (1-8m)}$, if $n$ is odd.
\end{lemma}

\begin{proof}
Again the solid torus boundary component of $\widehat{W}_n$ is tiled
by squares of side length $\sqrt{2}$, by Lemma \ref{lemma:k_shape},
which we view with sides parallel to the $x$ and $y$ axes in the
Euclidean plane.  The slope $\bar{\mu}_n$ still runs once along a
2--punctured disk, which came from $W_1$ or $\overline{W}_1$, hence
runs from $(0,0)$ to $(2\sqrt{2}, -2\sqrt{2})$, on the Euclidean
plane, also by Lemma \ref{lemma:k_shape}.

First suppose $n$ is even, $n=2k$.  The slope $\bar{\lambda}_n$ will
be formed by stepping $k+1$ times horizontally (following the meridian
of $K$ in the cusp tiling), and stepping $k-1$ times vertically
(following the meridian of $\overline{K}$).  Hence it runs from
$(0,0)$ to $(\sqrt{2}(k+1), \sqrt{2}(k-1))$.

Thus in the case $n=2k$, $\bar{\lambda}_n + m\,\bar{\mu}_n$ runs from
$(0,0)$ to $(\sqrt{2}(k+1 + 2m), \sqrt{2}(k-1-2m))$, so has length
$$ \left( 2( k+1+2m)^2 + 2(k-1-2m)^2 \right)^{1/2} = \sqrt{n^2 +
  4(1+2m)^2}.
$$

When $n$ is odd, $n=2k+1$, the slope $\bar{\lambda}_n$ is formed by
stepping $k$ times horizontally and $k+1$ times vertically, hence runs
from $(0,0)$ to $(\sqrt{2}k, \sqrt{2}(k+1))$.  So the slope
$\bar{\lambda}_n + m\,\bar{\mu}_n$ runs from $(0,0)$ to $(\sqrt{2}(k +
2m), \sqrt{2}(k+1-2m))$, and thus has length
$$
\left( 2 (k+2m)^2 + 2(k+1-2m)^2 \right)^{1/2} = \sqrt{n^2 + 16m^2 + (1-8m)}.
$$
\end{proof}

\begin{theorem}\label{thm:vol-withtwist}
For $n\geq 7$, or $n\geq 5$ and $|m|\geq 1$, the volume of the
complement of the $n$--chain link with $r$ (signed) half--twists is at
least:

\begin{tabular}{ccl}
  $\displaystyle{n\, v_8\, \left(1 - \frac{4\pi^2}{ n^2 + 16m^2} \right)^{3/2} }$
  & \hspace{.2in} & if $n$ is even and $r=2m$ is even, \\

  $\displaystyle{n\, v_8\, \left(1 - \frac{4\pi^2}{n^2 + 16m^2 + 16m +4} \right)^{3/2}}$
  & & if $n$ is even and $r=2m+1$ is odd, \\

  $\displaystyle{n\, v_8\, \left(1 - \frac{4\pi^2}{n^2 + 16m^2
  +(1+8m)} \right)^{3/2}}$ & 
  & if $n$ is odd and $r=2m$ is even, \\

  $\displaystyle{n\, v_8\, \left(1 - \frac{4\pi^2}{n^2 + 16m^2 + (1-8m)} \right)^{3/2}}$ &
  & if $n$ is odd and $r=2m+1$ is odd.
\end{tabular}

\smallskip

In all cases, the volume of the complement of that $n$--chain link is
larger than the volume of $W_{n-1}$ whenever $n\geq 60$ or $|m| \geq
8$.

When $n$ lies between $5$ and $10$, inclusively, the volume of the
complement of that $n$--chain link is larger than the volume of
$W_{n-1}$ whenever $|m| \geq 5$.
\end{theorem}

\begin{proof}
The volume estimates come from combining Lemmas
\ref{lemma:slope-length-nohalftwist} and
\ref{lemma:slope-length-halftwist} with Theorem \ref{thm:fkp}, using
the fact that the volumes of $\widehat{W}_n$ and $\overline{W}_n$ are
both $n\,v_8$, as they are both obtained by gluing $n$ copies of
manifolds isometric to the Whitehead link complement along totally
geodesic 3--punctured spheres.  

Note that in all cases, the bound on the volume is minimized for the
integer $m$ when $m=0$, in which case the argument of Theorem
\ref{thm:not-minvolume} still shows that the volume is greater than that of
$W_{n-1}$ for $n\geq 60$.

Let $\ell(n,m)$ denote the length of the Dehn filling slope, that is,
$\ell(n,m)$ is one of the four functions of $(n,m)$ in
Lemmas \ref{lemma:slope-length-nohalftwist}
and \ref{lemma:slope-length-halftwist}.  The volume of the chain link
is guaranteed to be strictly greater than that of $W_{n-1}$ whenever
we have
$$n\,v_8 \left(1 - \frac{4\pi^2}{\ell(n,m)^2}\right)^{3/2} \geq (n-1)\,v_8,$$
or whenever the function
$$f(n,m) = \frac{n}{n-1}\left(1
- \frac{4\pi^2}{\ell(n,m)^2}\right)^{3/2} -1$$ is strictly greater
than $0$.  Notice that in all four cases for $\ell(n,m)$, the function
$f$ is increasing with $|m|$.  Hence for fixed $n$, to find where this
function is greater than $0$, it suffices to set the function equal to
zero and solve for $m$.

We do so, and after a calculation, find that the zeros for $m$ can all
be computed in terms of the function
$$R(n) = \frac{1}{2}\sqrt{\frac{\pi^2}{1-\left(\frac{n-1}{n}\right)^{2/3}} -
\frac{n^2}{4}}.
$$
That is, the zeros of $f$ for $m$ are given by:

\begin{tabular}{ccl}
  $\pm {R(n)}$ & \hspace{.2in} & if $n$ and $r=2m$ are even,\\
  $-1/2\pm {R(n)}$ & & if $n$ is even and $r=2m+1$ is odd,\\
  $-1/4 \pm {R(n)}$ & & if $n$ is odd and $r=2m$ is even,\\
  $1/4 \pm {R(n)}$ & & if $n$ is odd and $r=2m+1$ is odd.
\end{tabular}

Now, check that $R(n)>0$ for any integer $n$ between $1$ and $59$,
inclusively, and that $R(n)$ is maximized at $n\approx 29.6104$, with
maximum value approximately $7.36$.  Thus the maximum and minimum
possible values for the zeros of $f$ lie strictly between $-8$ and
$8$.

As for the $n$--chain links with $n$ between $5$ and $10$,
inclusively, substituting the particular value of $n$ ($n=5, 6, \dots,
10$) into $R(n)$ and finding the zeros of $f$ for $m$, we see that all
such zeros lie strictly between $-6$ and $6$.
\end{proof}

In fact, for particular values of $n$ between $5$ and $59$,
inclusively, one can check that often the zeros of $m$ lie in a more
restrictive region than between $-8$ and $8$.  However, the power of
Theorem \ref{thm:vol-withtwist} is that it reduces the problem of
determining whether any $n$--chain link may have volume smaller than
that of $W_{n-1}$ to the task of checking only finitely many examples.

As in Section \ref{sec:compute}, we may now check the volumes of only
finitely many examples using the computer.  These finitely many
examples are checked by performing Dehn filling along slopes of the
form $(1,m)$, $m=0, \pm 1, \pm 2, \dots, \pm 7$, on $98$ initial
manifolds, namely the manifolds $\widehat{W}_n$ and $\overline{W}_n$
for $n = 11, 12, \dots, 59$.  In fact, when $n$ is odd,
$\widehat{W}_n$ and $\overline{W}_n$ are isometric, by an orientation
reversing isometry, and so it suffices to check volumes of Dehn
fillings of $\widehat{W}_n$ alone in this case.  In the even case,
$\widehat{W}_n$ and $\overline{W}_n$ must be checked separately.
However, in this case the Dehn fillings of $\widehat{W}_n$ along
slopes $(1,m)$ and $(1,-m)$ are isometric, by an orientation reversing
isometry, and Dehn fillings of $\overline{W}_n$ along $(1,m)$ and $(1,
-(m-1))$ are isometric, by an orientation reversing isometry.
Finally, since $(1,0)$ Dehn filling on $\widehat{W}_n$ gives the
minimally twisted $n$--chain link, which was examined in the previous
section, we omit that case.  In total, this leaves $14$ Dehn fillings
to check on the $25$ initial manifolds $\widehat{W}_n$ for $n$ odd,
$7$ Dehn fillings to check on the $24$ initial manifolds
$\widehat{W}_n$ for $n$ even, and $8$ Dehn fillings to check on the
$24$ initial manifolds $\overline{W}_n$ for $n$ even, or $710$ volumes
to compute by computer.

We automated the process of computing volumes as follows.
Triangulations for initial manifolds $\widehat{W}_n$ and
$\overline{W}_n$ were generated by Schleimer using the program
Twister \cite{twister}, which computes SnapPea triangulations
for manifolds described from the viewpoint of the mapping class group.

We then ran the triangulations through Snap \cite{goodman:snap}, to
find volumes of the manifolds under appropriate Dehn fillings.  These
volumes were compared with those of $W_{n-1}$.  In all cases, the
volume of $W_{n-1}$ was strictly smaller.  The data generated is shown
in Tables \ref{table:dehn-filling-voldata-1},
\ref{table:dehn-filling-voldata-2}, and
\ref{table:dehn-filling-voldata-tw}. 

Again in order to convert these results into a rigorous proof, the
algorithms of Moser \cite{moser} and Milley \cite{milley:minvol} could
be applied directly.  However, as the triangulations of general
$n$--chain link complements are of similar complexity as those of the
minimally twisted chain links, which were too complex for the computer
implementation of these algorithms for $n$ larger than 25, we omitted
this step.

\begin{table}
  {\footnotesize

\begin{tabular}{ | c || c | c | c | c | c | c | c || c |}
\hline

$n$ & -7 & -6 & -5 & -4 & -3 & -2 & -1 & $\vol(W_{n-1})$ \\

\hline

11 & 39.792 & 39.635 & 39.405 & 39.057 & 38.532 & 37.780 & 36.924 & 36.639\\
12 & 43.456 & 43.309 & 43.098 & 42.790 & 42.340 & 41.714 & 40.991 & 40.302\\
13 & 47.059 & 46.897 & 46.666 & 46.338 & 45.879 & 45.291 & 44.715 & 43.966\\
14 & 50.731 & 50.580 & 50.371 & 50.080 & 49.683 & 49.182 & 48.670 & 47.630\\
15 & 54.338 & 54.175 & 53.952 & 53.651 & 53.260 & 52.802 & 52.402 & 51.294\\
16 & 58.015 & 57.865 & 57.663 & 57.396 & 57.054 & 56.654 & 56.284 & 54.958\\
17 & 61.628 & 61.468 & 61.258 & 60.988 & 60.657 & 60.300 & 60.014 & 58.622\\
18 & 65.309 & 65.162 & 64.972 & 64.731 & 64.439 & 64.120 & 63.847 & 62.286\\
19 & 68.927 & 68.773 & 68.579 & 68.340 & 68.063 & 67.782 & 67.571 & 65.950\\
20 & 72.612 & 72.471 & 72.294 & 72.079 & 71.831 & 71.575 & 71.369 & 69.613\\
21 & 76.234 & 76.089 & 75.912 & 75.701 & 75.469 & 75.246 & 75.087 & 73.277\\
22 & 79.922 & 79.788 & 79.626 & 79.436 & 79.225 & 79.019 & 78.859 & 76.941\\
23 & 83.549 & 83.413 & 83.252 & 83.068 & 82.873 & 82.695 & 82.572 & 80.605\\
24 & 87.238 & 87.113 & 86.965 & 86.797 & 86.618 & 86.450 & 86.325 & 84.269\\
25 & 90.869 & 90.744 & 90.598 & 90.438 & 90.274 & 90.129 & 90.033 & 87.933\\
26 & 94.559 & 94.443 & 94.309 & 94.161 & 94.009 & 93.871 & 93.771 & 91.597\\
27 & 98.195 & 98.079 & 97.949 & 97.809 & 97.670 & 97.551 & 97.474 & 95.260\\
28 & 101.885 & 101.777 & 101.657 & 101.527 & 101.397 & 101.282 & 101.202 & 98.924\\
29 & 105.524 & 105.418 & 105.301 & 105.179 & 105.061 & 104.963 & 104.900 & 102.588\\
30 & 109.214 & 109.115 & 109.006 & 108.892 & 108.781 & 108.685 & 108.619 & 106.252\\
31 & 112.857 & 112.760 & 112.655 & 112.549 & 112.448 & 112.366 & 112.314 & 109.916\\
32 & 116.546 & 116.455 & 116.358 & 116.257 & 116.162 & 116.081 & 116.025 & 113.580\\
33 & 120.192 & 120.104 & 120.010 & 119.916 & 119.830 & 119.760 & 119.718 & 117.244\\
34 & 123.881 & 123.797 & 123.710 & 123.621 & 123.538 & 123.469 & 123.423 & 120.907\\
35 & 127.529 & 127.448 & 127.365 & 127.282 & 127.208 & 127.149 & 127.113 & 124.571\\
36 & 131.217 & 131.141 & 131.062 & 130.984 & 130.911 & 130.852 & 130.813 & 128.235\\
37 & 134.868 & 134.794 & 134.719 & 134.646 & 134.582 & 134.531 & 134.500 & 131.899\\
38 & 138.554 & 138.485 & 138.414 & 138.344 & 138.281 & 138.230 & 138.197 & 135.563\\
39 & 142.207 & 142.140 & 142.073 & 142.009 & 141.952 & 141.908 & 141.882 & 139.227\\
40 & 145.893 & 145.829 & 145.765 & 145.704 & 145.648 & 145.604 & 145.575 & 142.891\\
41 & 149.547 & 149.486 & 149.426 & 149.369 & 149.319 & 149.281 & 149.259 & 146.554\\
42 & 153.232 & 153.174 & 153.116 & 153.061 & 153.013 & 152.974 & 152.949 & 150.218\\
43 & 156.888 & 156.833 & 156.778 & 156.727 & 156.684 & 156.650 & 156.631 & 153.882\\
44 & 160.572 & 160.519 & 160.466 & 160.417 & 160.374 & 160.340 & 160.319 & 157.546\\
45 & 164.229 & 164.179 & 164.129 & 164.084 & 164.045 & 164.016 & 163.999 & 161.210\\
46 & 167.912 & 167.863 & 167.816 & 167.772 & 167.734 & 167.704 & 167.685 & 164.874\\
47 & 171.570 & 171.524 & 171.480 & 171.439 & 171.405 & 171.379 & 171.364 & 168.538\\
48 & 175.252 & 175.208 & 175.165 & 175.125 & 175.091 & 175.064 & 175.048 & 172.202\\
49 & 178.912 & 178.869 & 178.829 & 178.792 & 178.762 & 178.739 & 178.726 & 175.865\\
50 & 182.592 & 182.552 & 182.512 & 182.477 & 182.446 & 182.422 & 182.408 & 179.529\\
51 & 186.253 & 186.214 & 186.177 & 186.144 & 186.117 & 186.096 & 186.085 & 183.193\\
52 & 189.933 & 189.895 & 189.859 & 189.827 & 189.800 & 189.778 & 189.765 & 186.857\\
53 & 193.593 & 193.558 & 193.525 & 193.495 & 193.470 & 193.452 & 193.441 & 190.521\\
54 & 197.273 & 197.238 & 197.206 & 197.176 & 197.151 & 197.133 & 197.121 & 194.185\\
55 & 200.934 & 200.901 & 200.871 & 200.844 & 200.822 & 200.805 & 200.796 & 197.849\\
56 & 204.612 & 204.581 & 204.551 & 204.524 & 204.502 & 204.485 & 204.474 & 201.512\\
57 & 208.275 & 208.245 & 208.217 & 208.192 & 208.172 & 208.157 & 208.149 & 205.176\\
58 & 211.952 & 211.923 & 211.896 & 211.871 & 211.851 & 211.836 & 211.826 & 208.840\\
59 & 215.615 & 215.587 & 215.561 & 215.539 & 215.521 & 215.507 & 215.500 & 212.504\\

\hline
\end{tabular}
}

  \smallskip
  \caption{Volumes of chain links obtained by Dehn filling
    $\widehat{W}_n$ along slope $s = (1,m)$, where $m$ is the integer
    at the top of the column, compared with $\vol(W_{n-1})$.}
  \label{table:dehn-filling-voldata-1}
\end{table}

\begin{table}
  {\footnotesize

\begin{tabular}{ | c || c | c | c | c | c | c | c || c |}
\hline

$n$ & 1 & 2 & 3 & 4 & 5 & 6 & 7 & $\vol(W_{n-1})$ \\

\hline

11 & 37.340 & 38.184 & 38.821 & 39.249 & 39.531 & 39.721 & 39.851 & 36.639\\
13 & 44.982 & 45.598 & 46.126 & 46.517 & 46.792 & 46.985 & 47.122 & 43.966\\
15 & 52.581 & 53.034 & 53.467 & 53.813 & 54.072 & 54.263 & 54.403 & 51.294\\
17 & 60.139 & 60.478 & 60.829 & 61.131 & 61.370 & 61.553 & 61.693 & 58.622\\
19 & 67.662 & 67.919 & 68.205 & 68.465 & 68.682 & 68.855 & 68.991 & 65.950\\
21 & 75.155 & 75.354 & 75.586 & 75.810 & 76.005 & 76.166 & 76.296 & 73.277\\
23 & 82.623 & 82.780 & 82.971 & 83.162 & 83.336 & 83.484 & 83.607 & 80.605\\
25 & 90.073 & 90.197 & 90.355 & 90.520 & 90.673 & 90.809 & 90.925 & 87.933\\
27 & 97.506 & 97.607 & 97.738 & 97.879 & 98.015 & 98.139 & 98.247 & 95.260\\
29 & 104.926 & 105.009 & 105.119 & 105.240 & 105.361 & 105.473 & 105.573 & 102.588\\
31 & 112.335 & 112.404 & 112.497 & 112.602 & 112.708 & 112.810 & 112.902 & 109.916\\
33 & 119.735 & 119.792 & 119.872 & 119.963 & 120.057 & 120.149 & 120.234 & 117.244\\
35 & 127.127 & 127.176 & 127.244 & 127.323 & 127.407 & 127.489 & 127.567 & 124.571\\
37 & 134.513 & 134.554 & 134.613 & 134.682 & 134.757 & 134.831 & 134.903 & 131.899\\
39 & 141.893 & 141.928 & 141.979 & 142.040 & 142.106 & 142.174 & 142.240 & 139.227\\
41 & 149.268 & 149.299 & 149.343 & 149.397 & 149.456 & 149.517 & 149.577 & 146.554\\
43 & 156.639 & 156.665 & 156.704 & 156.752 & 156.805 & 156.860 & 156.916 & 153.882\\
45 & 164.006 & 164.029 & 164.064 & 164.106 & 164.154 & 164.204 & 164.254 & 161.210\\
47 & 171.370 & 171.390 & 171.421 & 171.459 & 171.501 & 171.547 & 171.594 & 168.538\\
49 & 178.731 & 178.749 & 178.776 & 178.810 & 178.849 & 178.890 & 178.933 & 175.865\\
51 & 186.089 & 186.106 & 186.130 & 186.160 & 186.195 & 186.233 & 186.272 & 183.193\\
53 & 193.446 & 193.460 & 193.482 & 193.509 & 193.541 & 193.576 & 193.612 & 190.521\\
55 & 200.800 & 200.813 & 200.832 & 200.857 & 200.886 & 200.918 & 200.951 & 197.849\\
57 & 208.152 & 208.164 & 208.181 & 208.204 & 208.230 & 208.259 & 208.290 & 205.176\\
59 & 215.503 & 215.513 & 215.529 & 215.550 & 215.574 & 215.601 & 215.629 & 212.504\\

\hline

\end{tabular}

}

  \smallskip
  \caption{Volumes of chain links obtained by Dehn filling
    $\widehat{W}_n$ along slope $s = (1,m)$, where $m$ is the integer
    at the top of the column, compared with $\vol(W_{n-1})$.}
  \label{table:dehn-filling-voldata-2}
\end{table}

\begin{table}
  {\footnotesize

\begin{tabular}{ | c || c | c | c | c | c | c | c |c || c |}

\hline

$n$ & -7 & -6 & -5 & -4 & -3 & -2 & -1 & 0 & $\vol(W_{n-1})$\\

\hline

12 & 43.513 & 43.389 & 43.214 & 42.959 & 42.586 & 42.049 & 41.349 & 40.709 & 40.302\\
14 & 50.790 & 50.661 & 50.484 & 50.237 & 49.896 & 49.443 & 48.914 & 48.492 & 47.630\\
16 & 58.076 & 57.945 & 57.771 & 57.539 & 57.234 & 56.858 & 56.456 & 56.165 & 54.958\\
18 & 65.370 & 65.241 & 65.073 & 64.858 & 64.591 & 64.279 & 63.971 & 63.763 & 62.286\\
20 & 72.671 & 72.545 & 72.387 & 72.191 & 71.958 & 71.701 & 71.461 & 71.308 & 69.613\\
22 & 79.979 & 79.858 & 79.711 & 79.534 & 79.332 & 79.119 & 78.930 & 78.814 & 76.941\\
24 & 87.292 & 87.178 & 87.042 & 86.883 & 86.708 & 86.531 & 86.380 & 86.290 & 84.269\\
26 & 94.610 & 94.503 & 94.378 & 94.237 & 94.085 & 93.937 & 93.815 & 93.744 & 91.597\\
28 & 101.933 & 101.833 & 101.718 & 101.592 & 101.461 & 101.337 & 101.237 & 101.180 & 98.924\\
30 & 109.259 & 109.166 & 109.062 & 108.950 & 108.836 & 108.730 & 108.647 & 108.601 & 106.252\\
32 & 116.588 & 116.502 & 116.407 & 116.307 & 116.208 & 116.119 & 116.049 & 116.011 & 113.580\\
34 & 123.919 & 123.840 & 123.754 & 123.665 & 123.579 & 123.501 & 123.443 & 123.411 & 120.907\\
36 & 131.253 & 131.179 & 131.101 & 131.022 & 130.946 & 130.880 & 130.830 & 130.803 & 128.235\\
38 & 138.588 & 138.520 & 138.449 & 138.379 & 138.312 & 138.254 & 138.211 & 138.188 & 135.563\\
40 & 145.924 & 145.861 & 145.797 & 145.734 & 145.675 & 145.625 & 145.588 & 145.568 & 142.891\\
42 & 153.260 & 153.203 & 153.145 & 153.088 & 153.036 & 152.992 & 152.960 & 152.943 & 150.218\\
44 & 160.598 & 160.545 & 160.492 & 160.441 & 160.395 & 160.356 & 160.328 & 160.313 & 157.546\\
46 & 167.936 & 167.888 & 167.839 & 167.793 & 167.752 & 167.717 & 167.693 & 167.680 & 164.874\\
48 & 175.274 & 175.230 & 175.186 & 175.144 & 175.107 & 175.076 & 175.055 & 175.043 & 172.202\\
50 & 182.613 & 182.572 & 182.532 & 182.494 & 182.461 & 182.433 & 182.414 & 182.404 & 179.529\\
52 & 189.951 & 189.914 & 189.877 & 189.843 & 189.813 & 189.788 & 189.771 & 189.762 & 186.857\\
54 & 197.290 & 197.255 & 197.222 & 197.190 & 197.163 & 197.141 & 197.126 & 197.118 & 194.185\\
56 & 204.629 & 204.596 & 204.566 & 204.537 & 204.512 & 204.493 & 204.479 & 204.471 & 201.512\\
58 & 211.967 & 211.937 & 211.909 & 211.883 & 211.860 & 211.843 & 211.830 & 211.824 & 208.840\\

\hline
\end{tabular}

}

  \smallskip
  \caption{Volumes of chain links obtained by Dehn filling
    $\overline{W}_n$ along slope $s = (1,m)$, where $m$ is the integer
    at the top of the column, compared with $\vol(W_{n-1})$.}
  \label{table:dehn-filling-voldata-tw}
\end{table}

\subsection{Chain links with 5 through 10 link components}

Our methods can be used to show that of all $n$--chain links, only the
minimally twisted $n$--chain link can possibly be the minimal volume
manifold with $n$ cusps for $n$ between $5$ and $10$, inclusively.

In fact, because the complexity of these manifolds was comparatively
small, we ran them through Milley's algorithm \cite{milley:minvol}, to
rigorously check this fact.  The algorithm successfully applied, and
we have the following theorem.

\begin{theorem}\label{thm:small-chains}
Let $n$ be an integer between $5$ and $10$, inclusively.  If $C_n$ is
an $n$--chain link that is \emph{not} minimally twisted, then the
complement $S^3\setminus C_n$ cannot be the minimal volume $n$--cusped
hyperbolic manifold.
\end{theorem}

\begin{proof}
Theorem \ref{thm:vol-withtwist} implies that for these $n$, and those
chain links with at least $11$ half--twists, the volume is strictly
greater than that of the $(n-1)$--fold cyclic cover over a component
of the Whitehead link, which is known to have larger volume than that
of the minimally twisted $n$--chain links in these cases.

The remaining cases to check are Dehn fillings $(1, \pm 1), (1, \pm
2), \dots, (1, \pm 5)$ on manifolds $\widehat{W}_n$ for $n$ odd,
Dehn fillings $(1,1), \dots, (1,5)$ on manifolds $\widehat{W}_n$ for
$n$ even, and Dehn fillings $(1,0), (1,1), \dots, (1,5)$ on manifolds
$\overline{W}_n$ for $n$ even.  These cases were run through
algorithms of Moser \cite{moser} and Milley \cite{milley:minvol}, and
their programs rigorously proved that the volumes of these chain links
were larger than that of the corresponding minimally twisted
$n$--chain link.  Programs are available from
Milley \cite{milley:minvol} or the second author.
\end{proof}

In Table \ref{table:small}, we show the volumes of $n$--chain link
complements, $n$ between $5$ and $10$, whose volumes are not
automatically larger than the minimally twisted $n$--chain link by
Theorem \ref{thm:vol-withtwist}.  These are compared with the volume
of the minimally twisted chain link.

\begin{table}
  {\footnotesize

\begin{tabular}{ | c ||  c | c | c | c | c |c || c |}

\hline

Fill mfld &  -5 & -4 & -3 & -2 & -1 & 0 & min twist \\
\hline

$\widehat{W}_5$ & 17.806 & 17.527 & 16.973 & 15.743
& 12.845 & - 
& 10.149\\ 
$\widehat{W}_6$ & 21.438 & 21.171 & 20.675 & 19.678
& 17.628 & - 
& 14.655\\ 
$\widehat{W}_7$ & 24.972 & 24.641 & 24.042 & 22.916
& 20.924 & - 
& 19.797\\ 
$\widehat{W}_8$ & 28.630 & 28.327 & 27.809 & 26.904
& 25.418 & - 
& 24.092\\ 
$\widehat{W}_9$ & 32.172 & 31.821 & 31.242 & 30.301
& 28.996 & - 
& 28.476\\ 
$\widehat{W}_{10}$ & 35.851 & 35.536 & 35.041 &
34.274 & 33.236 & - 
& 32.552\\ 
$\overline{W}_6$ &  21.526 & 21.324 & 20.963 & 20.266 & 18.832 & 16.000 & 14.655\\
$\overline{W}_8$  & 28.734 & 28.498 & 28.104 & 27.419 & 26.237 & 24.553 & 24.092\\
$\overline{W}_{10}$ & 35.963 & 35.711 & 35.317 & 34.696 & 33.776 & 32.759 & 32.552\\
\hline
\end{tabular}

\smallskip

\begin{tabular}{ | c ||  c | c | c | c | c || c |}

\hline

Fill mfld & 1 & 2 & 3 & 4 & 5 &  min twist  \\
\hline

$\widehat{W}_5$ & 14.603 & 16.485 & 17.301 & 17.688 & 17.894 & 10.149 \\
$\widehat{W}_6$ & 17.628 & 19.678 & 20.675 & 21.171 & 21.438 & 14.655 \\
$\widehat{W}_7$ & 22.040 & 23.567 & 24.387 & 24.830 & 25.082 & 19.797 \\
$\widehat{W}_8$ & 25.418 & 26.904 & 27.809 & 28.327 & 28.630 & 24.092 \\
$\widehat{W}_9$ & 29.672 & 30.824 & 31.568 & 32.018 & 32.294 & 28.476 \\
$\widehat{W}_{10}$ & 33.236 & 34.274 & 35.041 & 35.536 & 35.851 & 32.552 \\
$\overline{W}_6$ & 16.000 & 18.832 & 20.266 & 20.963 & 21.324 & 14.655 \\
$\overline{W}_8$ & 24.553 & 26.237 & 27.419 & 28.104 & 28.498 & 24.092 \\
$\overline{W}_{10}$ & 32.759 & 33.776 & 34.696 & 35.317 & 35.711 & 32.552 \\

\hline
\end{tabular}

}

  \smallskip
  
  \caption{Volumes of small chain links obtained by Dehn filling
    $\widehat{W}_n$ or $\overline{W}_n$ along slope $s = (1,m)$, where
    $m$ is the integer at the top of the column, compared with the
    volume of the minimally twisted chain
    link.}

  \label{table:small}
\end{table}

\bibliographystyle{amsplain}
\bibliography{../references}

\end{document}